\documentclass[a4paper,11pt]{amsart}

\usepackage{amssymb,mathrsfs,bm}
\usepackage{mathtools,galois,braket,esint}
\usepackage[utf8]{inputenc}
\usepackage[margin=2cm]{geometry}
\usepackage{hyperref,hypcap}
\hypersetup{colorlinks=true,allcolors=black,bookmarksdepth=3}
\usepackage[inline]{enumitem}
\usepackage{array}
\usepackage[utf8]{inputenc}
\usepackage{tikz, tikz-3dplot}
\usepackage{amsmath, amsthm, amssymb,graphics }
\usetikzlibrary{arrows,shapes,positioning}
\usetikzlibrary{intersections}
\usetikzlibrary{decorations.pathreplacing,calligraphy}

\usepackage{amsmath,amsthm,amssymb,amscd,color, xcolor,mathtools,url,tikz}
\usepackage{bbm}
\usepackage{pgfplots}

\theoremstyle{plain}
\newtheorem{theorem}{Theorem}[section]

\theoremstyle{definition}

\newtheorem{conjecture}[theorem]{Conjecture}
\theoremstyle{remark}

\numberwithin{equation}{section}
\numberwithin{figure}{section}
\numberwithin{table}{section}
\newcommand{\diff}{\mathop{}\!\mathrm{d}}

\newcommand{\Exp}{\mathrm{Exp}}

\newcommand{\dist}{\mathrm{dist}}
\newcommand{\vol}{\mathrm{vol}}

\DeclareMathOperator{\card}{card}

\allowdisplaybreaks

\begin{document}

\title{On universality for the kinetic wave equation}

\author[P.~Germain]{Pierre Germain}
\address[Pierre Germain]{Imperial College London Department of Mathematics}
\email{p.germain@imperial.ac.uk}
\author[H.~Zhu]{Hui Zhu}
\address[Hui Zhu]{Imperial College London Department of Mathematics}
\email{hui.zhu@imperial.ac.uk}

\begin{abstract}
    On compact Riemannian manifolds with chaotic geometries, specifically those exhibiting the random wave model conjectured by Berry, we derive heuristically a homogeneous kinetic wave equation that is universal for all such manifolds.
\end{abstract}

\maketitle


\section{Introduction}

\subsection{The kinetic wave equation}

Classical references on the theory of the \emph{kinetic wave equation} and its development can be found in the textbooks~\cite{Galtier, Nazarenko, ZLF}.
Although this theory is highly versatile and is applicable to any nonlinear dispersive equation, we shall focus on the model case of the \emph{nonlinear Schr\"odinger equation} posed on a Riemannian manifold $M$:
\begin{equation}
    \tag{NLS}
    \label{eq::NLS}
    i \partial_t u + \Delta u = \varepsilon |u|^2 u.
\end{equation}

Our attention is directed towards the \emph{weakly nonlinear setting}, where the strength of nonlinearity, denoted as $\varepsilon \in \mathbb{R} \setminus \{0\}$, is small. 
However, the sign of $\varepsilon$ holds no significance in our analysis. 
We are specifically interested in the long-time dynamics of the correlation function $\mathbb{E}[u(t,x)\overline{u(t,y)}]$, with $u$ representing a weakly interacting \emph{random wave ensemble}.
Assuming the \emph{phase randomization condition} of the waves, it has been conjectured that the dynamics obey a certain kinetic wave equation at the \emph{kinetic time order} (or the Van Hove limit \cite{VanHove1955}) $t \sim \varepsilon^{-2}$.

Take the torus $\mathbb{T}^d_L = \mathbb{R}^d / L \mathbb{Z}^d$, where $L>0$, for example.
The phase randomization condition mandates the Fourier coefficients
\begin{equation*}
    \hat{u}(t,k) = \frac{1}{L^{d/2}} \int e^{-2\pi i k \cdot x} u(t,x) \diff x
\end{equation*}
to be \emph{uncorrelated}.
For initially \emph{prepared} solutions --- these are solutions of the form $\hat{u}(0,k)=\phi(k)e^{i\theta_k}$ where $\phi$ decays sufficiently fast and $e^{i\theta_k}$ are independent and i.i.d.\ random phases --- this condition is proven to propagate in time (see e.g., \cite{DengHani2021propagation}).
In this context, the correlation function is entirely determined by the \emph{energy spectra} $\mathbb{E}|\hat{u}(t,k)|^2$.

Let us introduce the \emph{slow time} variable $\tau = t/T_{\mathrm{kin}}$ with
\begin{equation}
    \label{eq::Tkin}
    T_{\mathrm{kin}} = \pi \varepsilon^{-2}
\end{equation}
representing the kinetic order.
It has been demonstrated (see e.g., \cite{DengHani0,DengHani1,CollotGermain1,CollotGermain2}) that, taking the \emph{infinite volume limit} $L \to \infty$ and the \emph{weak nonlinearity limit} $\varepsilon \to 0$ simultaneously while adhering to certain scaling laws results in $\mathbb{E}|\hat{u}(\tau T_{\mathrm{kin}},k)|^2 \approx f(\tau,k)$, where $f$ satisfies the kinetic wave equation:
\begin{equation}
    \tag{KW}
    \label{eq::KW}
    \partial_\tau f(\tau,k) = \mathcal{C} [f](\tau,k),
\end{equation}
with initial data $f(0,k) = |\phi(k)|^2$.
The collision operator $\mathcal{C} [f]$ is expressed as
\begin{align*}
    \mathcal{C}[f](\tau,k) & = \iiint_{(\mathbb{R}^d)^3}  \mathcal{K}(k,k_1,k_2,k_3) \sum_{0\le j \le 3}(-1)^j \prod_{0 \le \ell \le 3}^{\ell \ne j} f(\tau,k_\ell) \diff k_1 \diff k_2 \diff k_3
\end{align*}
where $k_0 = k$, and the kernel $\mathcal{K}$ is given by
\begin{equation}
    \label{eq::KW-kernel}
    \mathcal{K}(k,k_1,k_2,k_3) = \bm{\delta}_{\mathbb{R}^d}(k - k_1 + k_2 - k_3) \bm{\delta}_{\mathbb{R}}(|k|^2 - |k_1|^2 + |k_2|^2 - |k_3|^2).
\end{equation}

The \emph{Kolmogorov--Zakharov (KZ) spectrum}, a stationary solution of the kinetic wave equation~\eqref{eq::KW}, holds particular significance as it is anticipated to depict turbulent cascades. Consequently, it often constitutes the primary focus of interest in various applications. Experimental validations, both physical and numerical, of equation~\eqref{eq::KW} predominantly rely on the observation of the KZ spectrum, see for instance reviews \cite{Galtier0,Nazarenko0,NewellRumpf} and recent developments \cite{BD,HrabskiPan,NEZLG,ZSKN}.

Remarkably, a series of mathematical works over recent years \cite{CollotGermain1,CollotGermain2,DengHani0,DengHani1,DengHani2,HaniShatahZhu,StaffilaniTran} culminated in a rigorous derivation of~\eqref{eq::KW}. This convergence of efforts has significantly advanced our understanding of the kinetic wave equation, providing a solid mathematical foundation for its theoretical underpinnings.

\subsection{The random wave model}

Known derivations, be they heuristic or rigorous, of the kinetic wave equation always take as their starting point Hamiltonian equations set in Euclidean geometry (full space or its quotients), the only exception being \cite{DGH}, which considers a model based on random matrices. This seems a considerable limitation on the validity of the kinetic wave equation; nevertheless, the KZ spectrum seems to have broad relevance, and it is natural to wonder whether an avatar of~\eqref{eq::KW} might hold in more general geometries.

We will focus in this article on \textit{chaotic geometries}; these are compact manifolds, with or without boundary (even though we will perform our analysis only on manifolds without boundaries), for which the geodesic flow is chaotic. Typically, \emph{chaotic billiards} and \emph{compact hyperbolic surfaces} fall under this category.
In this setting, Berry \cite{Berry1977,Berry1983} (see also Voros \cite{Voros2005}) conjectured that, in the semiclassical limit $\lambda \to \infty$, a corresponding eigenfunction \emph{$\psi$ ``appears to be a (real) Gaussian random function exhibiting ... oscillations ... (that) are statistically isotropic with the autocorrelation ... given by a Bessel function.''}.
This is the \emph{random wave model (RWM)}.

To date, different mathematical formulations of Berry's conjecture exist, see e.g., \cite{RS1994hyperbolic,Zelditch2010recent,NazarovSodin2010random,Nonnenmacher2013anatomy,Bourgain,Abert2018,Ingremeau,Sodin2016lectures}.
To illustrate the ideas, let us follow Bourgain \cite{Bourgain} (see also \cite{Ingremeau}) who states that, in the semiclassical limit $\lambda \to \infty$, the ensemble of functions $(\psi^x)_{x\in M}$, where $\psi^x$ is the restriction of $\psi$ to a geodesic ball centered at $x$ with radius much larger than the Planck scale $1/\sqrt{\lambda}$, resembles the \emph{monochromatic Gaussian random wave} $\Phi$ on $\mathbb{R}^d$.
Recall that, this $\Phi$ is an isotropic Gaussian random field on $\mathbb{R}^d$ whose spectrum is given by
\begin{equation*}
    \frac{1}{(2\pi)^d} \int e^{-i \xi \cdot Z} \mathbb{E}\Bigl[\Phi\Bigl(X+\frac{Z}{2}\Bigr) \Phi\Bigl(X-\frac{Z}{2}\Bigr)\Bigr] \diff Z
    = \frac{\bm{\delta}_{\mathbb{S}^{d-1}}(\xi)}{s(d)}
\end{equation*}
where $s(d) = 2\pi^{d/2} \Gamma(d/2)^{-1}$ is the area of $\mathbb{S}^{d-1}$.
Equivalently, this means that the correlation function of $\Phi$ is a radial function $\Lambda_d$ given by the Fourier transform
\begin{equation}
    \label{eq:def-J}
    \Lambda_d\bigl(q\bigr) = \frac{1}{s(d)} \int_{\mathbb{S}^{d-1}} e^{iZ \cdot \theta} \diff \theta,
\end{equation}
where $Z \in \mathbb{R}^d$ is any vector with $|Z|=q$.
It is know that $\Lambda_d$ can be computable in terms of Bessel functions (see e.g., \cite[Ch.~\textrm{VIII}, Eq.~(25)]{Stein}):
$\Lambda_d(q) = \Gamma(d/2) |q/2|^{-\nu} J_{\nu}(q)$,
where $ \nu = d/2-1$ and $J_\nu$ is the Bessel function of the first kind of order $\nu$.

Following Berry, to make sense of this approximation, one needs to define the (multi-point) correlations of the ensemble by taking local averages over domains much larger than the Planck scale $1/\sqrt{\lambda}$.
We postpone till \S\ref{sec::RWM} a discussion of the properties of the RWM in greater detail, but immediately provide some references pointing towards its validity. Rigorous results on the RWM are scant \cite{Bourgain,Ingremeau}, but physical
\cite{HSTBS,SHS} and numerical \cite{HR,MK} evidence is available, in particular since the connection of the RWM to the number of nodal domains was established \cite{BS}.

In the present article, we propose a heuristic derivation of a kinetic wave equation, which is valid under Berry's conjecture.

\subsection{Main result: kinetic limit on domains with chaotic geometries}

\label{sec::main-result}

We will now state precisely our main (formal) result: the solutions of~\eqref{eq::NLS} on a domain with chaotic geometry converge to a universal kinetic wave equation.

Consider a compact manifold $(M,g)$ with chaotic geometry and without boundary, for example a compact hyperbolic surface.
Let $0 = \lambda_0 \le \lambda_1 \le \lambda_2 \le \cdots$ be eigenvalues (counting multiplicity) of $-\Delta$ and let $\psi_n$ be the corresponding eigenfunctions that are orthogonal.
We further assume that the eigenfunctions are real-valued (even though this condition is not necessary for our results) and are normalized with $\|\psi_n\|_{L^2(M)}^2 = \vol(M)$.
For $L > 0$, on the dilated manifold $M_L = (M,L^2 g)$, one has eigenvalues and eigenfunctions $(\lambda_n^L,\psi_n)$, with $\lambda_n^L = \lambda_n / L^2$.

We set~\eqref{eq::NLS} on $M_L$, and provide it with prepared initial data
\begin{equation}
    \label{eq::ini-data-NLS}
    u(0,x) = \frac{1}{L^{d/2}} \sum_{n \ge 0} \varphi(\lambda_n^L) e^{i\vartheta_n} \frac{\psi_n(x)}{\sqrt{\vol(M)}},
\end{equation}
where $\vartheta_n$ are i.i.d.\ random variables, uniformly distributed on $[0,2\pi]$, and $\varphi$ is a smooth, rapidly decaying function.
Our normalization corresponds to $\| u(0,\cdot) \|_{L^2(M_L)}^2 =\mathcal{O}(1)$ and $\| u(0,\cdot) \|_{L^\infty}= \mathcal{O}(L^{-d})$ up to logarithmic corrections.
The solution $u$ to~\eqref{eq::NLS} can be expanded in the orthonormal basis as
\begin{equation}
    \label{eq::sol-NLS-fourier}
    u(t,x) = \frac{1}{L^{d/2}} \sum_{n \ge 0} A_n(t) \frac{\psi_n(x)}{\sqrt{\vol(M)}}.
\end{equation}
Under the simultaneous infinite volume limit $L \to \infty$ and the weak nonlinearity limit $\varepsilon \to 0$, we assert that $\mathbb{E} |A_n|^2 (\tau T_{\mathrm{kin}}) \approx \rho(\tau,\omega)$ as $\lambda_n^L \to \omega$.
Here $\rho$ solves the kinetic wave equation
\begin{equation}
    \label{eq::KWR}
    \tag{KWR}
    \partial_\tau \rho(\tau,\omega) = \mathcal{C}_* [\rho](\tau,\omega),
\end{equation}
with initial data $\rho(0,\omega) = |\varphi(\omega)|^2$.
Identifying $\omega = \omega_0$, the collision operator $\mathcal C_*[\rho]$ is expressed as
\begin{align*}
    \mathcal{C}_*[\rho](\tau,\omega) & =  \iiint_{(\mathbb{R}_+)^3} \mathcal{K}_*(\omega, \omega_1 , \omega_2 , \omega_3)
    \sum_{0\le j \le 3}(-1)^j \prod_{0 \le \ell \le 3}^{\ell \ne j} \rho(\tau,\omega_\ell) \diff \omega_1 \diff \omega_2 \diff \omega_3
\end{align*}
where the kernel $\mathcal{K}_*$ is defined by
\begin{equation*}
    \mathcal{K}_*(\omega,\omega_1,\omega_2,\omega_3)
    = \frac{\pi^2}{2} \biggl(\frac{s(d)}{(2\pi)^d}\biggr)^3 (\omega_1\omega_2\omega_3)^{d/2-1} \mathcal{I}(\omega,\omega_1,\omega_2,\omega_3) \bm{\delta}_{\mathbb{R}}(\omega - \omega_1 + \omega_2 - \omega_3),
\end{equation*}
and the interaction integral $\mathcal{I}$ is defined by
\begin{equation*}
    \mathcal{I}(\omega,\omega_1,\omega_2,\omega_3)
    = s(d) \int_{\mathbb{R}_+} q^{d-1} \prod_{0 \le j \le 3} \Lambda_d(\sqrt{\omega_j} q) \diff q.
\end{equation*}
As we shall see in \S\ref{sec::interaction-integral}, this integral $\mathcal{I}(\omega,\omega_1,\omega_2,\omega_3)$ quantifies the interactions among four eigenfunctions, serving the same role as $\delta_{\mathbb{R}^d}(k-k_1+k_2-k_3)$ in~\eqref{eq::KW-kernel}.

Detailed computations leading to the derivation of~\eqref{eq::KWR} are provided in \S\ref{sec::derivation}.
We shall also see in \S\ref{sec::KW-to-KWR} that the equation~\eqref{eq::KWR} also describes the dynamics of radial solutions to~\eqref{eq::KW}.
Therefore, we have named the kinetic equation~\eqref{eq::KWR} due to the dual interpretation of the letter \emph{``R''} not only as the \emph{``random wave model''} but also as \emph{``radial''}.

\section{Discussion of the result}
\subsection{Range of validity of the singular limit}

The derivation of \eqref{eq::KWR} in \S\ref{sec::derivation} relies on the convergence of the sum over frequencies to the integral defining the collision operator, see \S\ref{sec::kinetic-limit} for the details.
To ensure the validity of this sum-to-integral limit, the points $(\lambda_n^L,\lambda_{n_1}^L, \lambda_{n_2}^L, \lambda_{n_3}^L)$ (with $n$ fixed) must be equidistributed in a $T_{\mathrm{kin}}^{-1}$ neighborhood around the resonant manifold defined by the equation $\Omega(\omega, \omega_1,\omega_2,\omega_3) = 0$ (with $\omega$ fixed).

Here, all frequencies $\lambda_n^L$ can be confined to be of size $\mathcal{O}(1)$ due to the rapid decay of $\varphi$. 
Furthermore, we can assume heuristically that the frequencies are random, uniformly distributed and independent.
Thus, it suffices for there to be an infinite number of points in the vicinity of the resonant manifold.
Considering the volume of this neighborhood, which is of order $T_{\mathrm{kin}}^{-1}$, and the number of points in a unit volume, which is of order $L^{3d}$, we obtain the condition $L^{3d} T_{kin}^{-1} \gg 1$, or equivalently $L^{-3d/2} \ll \varepsilon$.

Another condition that is necessary to ensure the convergence of the $\operatorname{sinc}^2$ function to the Dirac $\bm{\delta}_{\mathbb{R}}$, see also \S\ref{sec::derivation}.
This convergence occurs as $T_{\mathrm{kin}} \to \infty$, i.e., when $\varepsilon \to 0$.

In summary, the kinetic limit is deemed valid in the limit where
\begin{equation*}
    L^{-3d/2} \ll \varepsilon \ll 1.    
\end{equation*}

\subsection{\texorpdfstring{From~\eqref{eq::KW} to~\eqref{eq::KWR}}{From (KW) to (KWR)}}

\label{sec::KW-to-KWR}

We observe that~\eqref{eq::KWR} can also be derived as the radial form of~\eqref{eq::KW}. Indeed, if $f$ is radial, choosing radial coordinates $2\pi k = \sqrt{\omega} \theta$ where $\omega \ge 0$, $\theta \in \mathbb{S}^{d-1}$, and setting $\rho(\tau,\omega) = f(\tau,k)$, we claim that
\begin{equation*}
    \mathcal{C}_*[\rho](\tau,\omega)
    = \mathcal{C}[f](\tau,k) = \frac{1}{s(d)} \int_{\mathbb{S}^{d-1}} \mathcal{C}[f]\bigl(\tau,|k|\theta\bigr) \diff \theta.
\end{equation*}
To show this, let us write
\begin{align*}
    \mathcal{C}[f](\tau,k)
     & = \iiint_{(\mathbb{R}^d)^3}  \mathcal{K}\Bigl(\frac{\sqrt{\omega}\theta}{2\pi},\frac{\sqrt{\omega_1}\theta_1}{2\pi},\frac{\sqrt{\omega_2}\theta_2}{2\pi},\frac{\sqrt{\omega_3}\theta_3}{2\pi}\Bigr) \sum_{0\le j \le 3}(-1)^j \prod_{0 \le \ell \le 3}^{\ell \ne j} \rho(\tau,\omega_\ell) \diff k_1 \diff k_2 \diff k_3.
\end{align*}
Using $\diff k = g(\omega) \diff \omega \diff \theta$ where $g(\omega) = \frac{1}{2} (2\pi)^{-d} \omega^{d/2-1}$ and expressing
\begin{equation*}
    \bm{\delta}_{\mathbb{R}^d}\Bigl(\frac{\sqrt{\omega}\theta}{2\pi}-\frac{\sqrt{\omega_1}\theta_1}{2\pi}+\frac{\sqrt{\omega_2}\theta_2}{2\pi}-\frac{\sqrt{\omega_3}\theta_3}{2\pi}\Bigr)
    = \int_{\mathbb{R}^d} e^{iz\cdot(\sqrt{\omega} \theta - \sqrt{\omega_1} \theta_1 + \sqrt{\omega_2} \theta_2 - \sqrt{\omega_3} \theta_3)}  \diff z,
\end{equation*}
our claim follows from a direct computation that shows:
\begin{align*}
    \mathcal{K}_* & (\omega,\omega_1,\omega_2,\omega_3)                                                                                                                                                                                                                                                                                 \\
                  & = \frac{g(\omega_1) g(\omega_2) g(\omega_3)}{s(d)} \iiiint_{(\mathbb{S}^{d-1})^4} \mathcal{K}\Bigl(\frac{\sqrt{\omega}\theta}{2\pi},\frac{\sqrt{\omega_1}\theta_1}{2\pi},\frac{\sqrt{\omega_2}\theta_2}{2\pi},\frac{\sqrt{\omega_3}\theta_3}{2\pi}\Bigr) \diff \theta \diff \theta_1 \diff \theta_2 \diff \theta_3.
\end{align*}

It is well-known that, in dimension 3 and for radial functions, the collision operator takes a very simple form, see for instance the appendix of~\cite{ZSKN}. By the above relation, it is also true for~\eqref{eq::KWR}.

\subsection{Universality of turbulent spectra}

We showed that the kinetic wave equation is universal for domains with chaotic geometry. Since the KZ spectra arise as stationary solutions of~\eqref{eq::KWR}, they are also universal over this class of domains.

Note that the ideas used in the present article are applicable to any isotropic equation, and thus that universality of turbulent spectra should extend to this class of equations.

\subsection{Questions}

To summarize the findings of the present article, the kinetic equation~\eqref{eq::KWR} is valid on domains with chaotic geometry, but also on tori under an  isotropy assumption (which we do not attempt to formulate more precisely).

It is tempting to ask if the class of domains for which~\eqref{eq::KWR} is valid might not be even larger.
Amongst Riemannian manifolds, the sphere often exhibits the most irregular behavior in terms of equidistribution of eigenvalues over $\mathbb{R}$, or equidistribution of eigenfunctions over the manifold. Therefore, the sphere seems to be the most natural example for which~\eqref{eq::KWR} might not apply.

One can also ask the same question beyond the framework of Riemannian manifolds: what can be said of confining potentials? Or even of graphs in the proper limit?

\section{The random wave model and its implications}

\label{sec::RWM}

In this section, we discuss the random wave model, review and mathematically formulate Berry's conjecture, and obtain several consequences of it.
We shall henceforth denote $f(n) \sim g(n)$ for two functions of variable $n \in \mathbb{N}$, when $g(n) \ne 0$ for $n \gg 1$ and $\lim_{n \to \infty} f(n) / g(n) = 1$.
We shall also denote $f(n) \sim 0$ when $\lim_{n \to \infty} f(n) = 0$.

\subsection{Berry's conjecture}

\begin{conjecture}[Berry]
    When $n \gg 1$, the ensemble of random variables $(\psi^x_n)_{x\in M}$ with respect to the probability measure on $M$:
    \begin{equation*}
        \diff\mu = \vol(B)^{-1} \bm{1}_{B_r} \diff x,
    \end{equation*}
    where $B_r$ is a geodesic ball of radius $r \gg 1/\sqrt{\lambda_n}$, resembles an isometric Gaussian random field on $\mathbb{R}^d$ with variance $\Lambda_d$, given by~\eqref{eq:def-J}.
\end{conjecture}

When the geometry is locally Euclidean (e.g., in the interior of a chaotic billiard in $\mathbb{R}^d$), this means that for any sequence $(z^j)_{1\le j\le m} $ of vectors in $\mathbb{R}^d$, the $m$-tuple of random variables $\psi_n^j$ defined by $\psi_n^j(x) = \psi_n(x+z^j) $ is approximately a joint Gaussian distribution with respect to $\diff \mu$.
This Gaussianity is characterized via the following Isserlis--Wick type formula for multi-point correlations:
\begin{equation}
    \label{eq::Berry-Wick}
    \Bigl\langle \prod_j \psi_n^j \Bigr\rangle_\mu
    \sim \sum_{P} \prod_{\{j,j'\}\in P} \langle \psi_n^j \psi_n^{j'} \rangle_\mu,
\end{equation}
where $\langle \cdot \rangle_\mu$ denote the expectation with respect to $\mu$, $P$ ranges among all partitions of $\{1,\ldots,m\}$ into disjoint doubletons, and the two-point correlations are given by
\begin{equation*}
    \langle \psi_n^j \psi_n^{j'} \rangle_\mu
    = \Lambda_d\bigl(\sqrt{\lambda_n}|z^{j}-z^{j'}|\bigr).
\end{equation*}

On a non-Euclidean manifold $M$ (which is compact and without boundary), we need $z^j:M \to TM$ to be (smooth) vector fields.
In this situation, the random variables $\psi_n^j$ are defined by $\psi_n^j(x) = (\psi_n \comp \Exp_x)(z^j_x)$, where $\Exp_x$ is the exponential map at $x$.
We conjecture that, the formula~\eqref{eq::Berry-Wick} still holds in this setting, with the two-point correlation given by
\begin{equation}
    \label{eq::Berry-two-point-M}
    \langle \psi_n^j \psi_n^{j'} \rangle_\mu
    = \bigl\langle \Lambda_d\bigl(\sqrt{\lambda_n}D_{j,j'}\bigr) \bigr\rangle_\mu, \quad D_{j,j'}(x) = \dist\bigl(\Exp_x(z_x^j),\Exp_x(z_x^{j'})\bigr).
\end{equation}
We believe that this formulation is essentially equivalent to the one given by Sodin in \cite{Sodin2016lectures}.
Particularly, when $m=2\ell$ is even and $z_j \equiv 0$, the estimate~\eqref{eq::Berry-Wick} yields
    \begin{equation*}
        \langle \psi^{2\ell} \rangle_{\mu}
        \sim \frac{\ell !}{2} \binom{2\ell}{\ell} = \frac{(2\ell)!}{2\ell!},
    \end{equation*}
the right hand side being the number of all possible ways of partitioning a set of $2\ell$ elements into disjoint doubletons.
Taking $\ell = 1$, one recovers the quantum unique ergodicity $\langle \psi^2 \rangle_{\mu} \sim 1$ at scale $r$.
We refer to works \cite{LR2017small,GI2017planck} for investigations of the quantum unique ergodicity at small scales on tori.

With the natural conjecture that $\psi_n$ and $\psi_m$ being uncorrelated with respect to $\mu$ when $n \ne m$, Berry's conjecture can be extended to the ensemble $(\psi_n^x)_{n\ge 0,x \in M}$.
In this case, under the assumptions that $n = \min n_j \gg 1$ and $r \gg 1/\sqrt{\lambda_n}$, the estimate~\eqref{eq::Berry-Wick} becomes
\begin{equation}
    \label{eq::Berry-Wick-M-multi-eigen}
    \Bigl\langle \prod_j (\psi_{n_j} \comp \Exp_x)(z_x^j) \Bigr\rangle_\mu
    \sim \sum_{P} \prod_{\{j,j'\}\in P} \bm{1}_{n_j=n_{j'}} \Bigl\langle \Lambda_d\Bigl(\sqrt{\lambda_{n_j}}D_{j,j'} \Bigr) \Bigr\rangle_\mu,
\end{equation}
where $D_{j,j'}$ is defined as in~\eqref{eq::Berry-two-point-M}.
Particularly, taking $z^j = 0$ in~\eqref{eq::Berry-Wick-M-multi-eigen} yields
\begin{equation}
    \label{eq::Berry-multi-eigen-z=0}
    \Bigl\langle \prod_j \psi_{n_j}\Bigr\rangle_\mu
    \sim \sum_{P} \prod_{\{j,j'\}\in P} \bm{1}_{n_j=n_{j'}}.
\end{equation}
As another consequence, if $m=2$, if $z:M\to SM$ and $q \in \mathbb{R}$, then~\eqref{eq::Berry-two-point-M} implies that
\begin{equation}
    \label{eq::Berry-two-point-M-Wigner}
    \Bigl\langle (\psi_n \comp \Exp_x)\Bigl(\frac{1}{2} q z_x\Bigr) \times (\psi_{n'} \comp \Exp_x)\Bigl(-\frac{1}{2} q z_x\Bigr)\Bigr\rangle_\mu \sim \bm{1}_{n=n'} \Lambda_d(\sqrt{\lambda_n} q).
\end{equation}

\subsection{Weak formulation}

For our derivation of the kinetic wave equation, as we shall see, a weak formula of Berry's conjecture is already sufficient.
We start from~\eqref{eq::Berry-two-point-M-Wigner} with the belief that local averages in $SM$ yields a similar estimate:
    \begin{equation}
        \label{eq::Berry-two-point-SM}
        \Bigl\langle (\psi_n \comp \Exp_x)\Bigl(\frac{1}{2} q \theta\Bigr) \times (\psi_{n'} \comp \Exp_x)\Bigl(-\frac{1}{2} q \theta\Bigr)\Bigr\rangle_{\nu} \sim \bm{1}_{n=n'}\Lambda_d(\sqrt{\lambda_n} q),
    \end{equation}
where $\diff\nu$ is a probability measure on $SM$ lifting $\diff \mu$:
    \begin{equation*}
        \diff\nu = \frac{\bm{1}_{(x,\theta) \in SB_r} \diff (x,\theta)}{\vol(SB_r)} \sim \frac{\bm{1}_{(x,\theta) \in SB_r} \diff (x,\theta)}{s(d) \vol(B_r)}.
    \end{equation*}
We can easily deduce~\eqref{eq::Berry-two-point-SM} from~\eqref{eq::Berry-two-point-M-Wigner} as follows.
Let $x_0$ be the center of the geodesic ball $B_r$.
For $\eta \in S_{x_0}M$, define the vector field $z^\eta : B_r \to TM$ by $z^\eta_x = \Gamma_x (\eta)$, where $\Gamma_x:T_{x_0}M \to T_xM$ is the parallel transport from $x_0$ to $x$ along the (radial) geodesic (with respect to the Levi--Civita connection).
Is is known that $\Gamma_x$ is a linear isometry, so $z^\eta_x \in S_xM$. Identifying $S_{x_0} \simeq \mathbb{S}^{d-1}$, we have
\begin{align*}
    \Bigl\langle (\psi_n & \comp \Exp_x)\Bigl(\frac{1}{2} q \theta\Bigr) \times (\psi_{n'} \comp \Exp_x)\Bigl(-\frac{1}{2} q \theta\Bigr)\Bigr\rangle_{\nu} \\
    & \sim \frac{1}{s(d)}\int_{\mathbb{S}^{d-1}} \Bigl\langle (\psi_n \comp \Exp_x)\Bigl(\frac{1}{2} q z^\eta_x\Bigr) \times (\psi_{n'} \comp \Exp_x)\Bigl(-\frac{1}{2} q z^\eta_x\Bigr)\Bigr\rangle_\mu \diff \eta
    \sim \bm{1}_{n=n'}\Lambda_d(\sqrt{\lambda_n} q).
\end{align*}
Similar analysis applies to multi-point correlations:
letting $(n_j)_{1\le j \le m}$ be distinct indices, so that no pairings in between them exist, then
\begin{equation}
    \label{eq::Berry-multi-point-SM}
    \Bigl\langle \prod_j (\psi_{n_j} \comp \Exp_x)\Bigl(\frac{1}{2} q \theta\Bigr) \times \prod_j (\psi_{n_j} \comp \Exp_x)\Bigl(-\frac{1}{2} q \theta\Bigr) \Bigr\rangle_{\nu}
    \sim \prod_j \Lambda_d\Bigl(\sqrt{\lambda_{n_j}} q\Bigr). 
\end{equation}

\subsection{Interaction integrals}

\label{sec::interaction-integral}

The principle term of the multi-point correlation~\eqref{eq::Berry-Wick-M-multi-eigen} vanishes when no pairings $P$ exists.
This happens particularly when $n_j$ are distinct.
However, we can still compute its sub-principal term, when $z_j \equiv 0$, using~\eqref{eq::Berry-multi-point-SM}.
In fact, using the spherical coordinates,
\begin{align}
    \label{eq:est-interaction-real-square}
    \vol(B_r) \Bigl\langle \prod_j \psi_{n_j} \Bigr\rangle_\mu^2
     & \sim s(d) \int_{\mathbb{R}^+} q^{d-1}
    \Bigl\langle \prod_j(\psi_{n_j} \comp \Exp_x)\Bigl(\frac{1}{2} q \theta\Bigr) \times \prod_j (\psi_{n_j} \comp \Exp_x)\Bigl(-\frac{1}{2} q \theta\Bigr) \Bigr\rangle_{\nu} \diff q \\
     & \sim s(d) \int_{\mathbb{R}^+} q^{d-1} \prod_j \Lambda_d\Bigl(\sqrt{\lambda_{n_j}} q\Bigr) \diff q
    \eqcolon \mathcal{I}(\lambda_{n_1},\ldots,\lambda_{n_m}).\nonumber
\end{align}
We shall call $\mathcal{I}(\lambda_{n_1},\ldots,\lambda_{n_m})$ an interaction integral.
It clearly has the scaling invariance:
\begin{equation}
    \label{eq::scaling-I}
    \mathcal{I}(\lambda_{n_1}/L^2,\ldots,\lambda_{n_m}/L^2)
    = L^d \mathcal{I}(\lambda_{n_1},\ldots,\lambda_{n_m}).
\end{equation}

One then expects estimates similar to~\eqref{eq:est-interaction-real-square} to hold over subsets $U \subset M$ much larger than the plank scale $1/\sqrt{\lambda}$, particularly when $U=M$.
Precisely, letting $\langle \cdot \rangle_U$ denote the expectation with respect to $\vol(U)^{-1} \bm{1}_{U} \diff x $, one expects that, when $n_j$ are distinct, the multi-point correlation $\langle \prod_j \psi_{n_j} \rangle_U$ is determined, up to a sign, by the square root of the corresponding interaction integral:
\begin{equation}
    \label{eq:est-interaction-real}
    \sqrt{\vol(U)} \Bigl\langle \prod_j \psi_{n_j} \Bigr\rangle_U
    \sim \pm \sqrt{\mathcal{I}(\lambda_{n_1},\ldots,\lambda_{n_m})}.
\end{equation}
To establish~\eqref{eq:est-interaction-real}, it is essential to ensure that distant points make negligible contributions to the correlation.
This is quantified in the following bound, which is suggested by the decay of $\Lambda_d$, but probably not optimal:
\begin{equation*}
    \iint \bm{1}_{\dist(x,y) > r} \prod_{1\le j \le m} \psi_{n_j}(x) \psi_{n_j}(y)  \diff x \diff y
    \lesssim \prod_j \Bigl(\sqrt{\lambda_{n_j}}r\Bigr)^{-(d-1)/2}.
\end{equation*}

\subsection{Phase randomization}

As stated by \cite[\S2.1.2]{ZLF}, due to the dispersion of the wave system, even initially correlated oscillations with different wavenumbers undergo phase randomization as time progresses.
Thus, in describing a free wave field it would be natural to average over the ensemble of chaotic (random) phases, i.e., to use the random phase approximation.
This assumption is closely related to the molecular chaos assumption in kinetic theory \cite{Boltzmann1872,Gallagher2013}.
Let us therefore choose i.i.d.\ random variables $\vartheta_n$, uniformly distributed on $[0,2\pi]$, and put $\Psi_n^\pm(x) = \psi_n(x) e^{\pm i\vartheta_n}$.
The correlation formulas can be stated in terms of $\Psi_n^\pm$.
The two-point correlation now states that
\begin{equation*}
    \mathbb{E}\Bigl\langle (\Psi_n^+ \comp \Exp_x)\Bigl(\frac{1}{2} q z_x\Bigr) \times (\Psi_{n'}^- \comp \Exp_x)\Bigl(-\frac{1}{2} q z_x\Bigr)\Bigr\rangle_\mu \sim \bm{1}_{n=n'} \Lambda_d\bigl(\sqrt{\lambda_n} q\bigr).
\end{equation*}
Similar analysis as in the previous subsection yields that: if $n_j$ are distinct, then
\begin{equation*}
    \vol(U)\, \mathbb{E}\Bigl| \Bigl\langle \prod_j \Psi_{n_j}^{\sigma_j} \Bigr\rangle_U\Bigr|^2
    = \vol(U)\, \Bigl| \Bigl\langle \prod_j \Psi_{n_j}^{\sigma_j} \Bigr\rangle_U\Bigr|^2
    \sim  \mathcal{I}(\lambda_{n_1},\ldots,\lambda_{n_m}).
\end{equation*}
Consequently, for some random phase $e^{i\vartheta}$ (in fact $\vartheta = \sum \sigma_j \vartheta_j$), there holds
\begin{equation*}
    \sqrt{\vol(U)}  \Bigl\langle \prod_j \Psi_{n_j}^{\sigma_j} \Bigr\rangle_U
    \sim \pm e^{i\vartheta} \sqrt{\mathcal{I}(\lambda_{n_1},\ldots,\lambda_{n_m})}.
\end{equation*}

\section{Derivation of the kinetic wave equation}

\label{sec::derivation}

In this section, we derive the kinetic wave equation~\eqref{eq::KWR}.
For this purpose, let us fix $(M,g)$ to be a compact manifold with chaotic geometry and without boundary, so that the analysis in \S\ref{sec::RWM} applies.

\subsection{Fourier series expansion}

Let the eigenvalues and eigenfunctions $(\lambda_n,\psi_n)$ be chosen as in \S\ref{sec::main-result}.
Recall that, by Weyl's law, we have $\card\{n:\lambda_n \le \lambda\} \sim (2\pi)^{-d} v(d) \vol(M) \lambda^{d/2}$, where $v(d) =\pi^{d/2} \Gamma(d/2+1)^{-1}$ is the volume of the unit ball in $\mathbb{R}^d$.
Therefore, setting $C_M = \vol(M)^{-1} \Gamma(d/2+1)$, one obtains the asymptotic behavior for eigenvalues:
\begin{equation}
    \label{eq::weyl-law}
    \lambda_n \sim 4\pi (C_M n)^{2/d}.
\end{equation}

Let $u$ be the solution to~\eqref{eq::NLS} on $M_L$ with initial data~\eqref{eq::ini-data-NLS}, and let $A_n$ be the Fourier coefficients of $u$ as determined by~\eqref{eq::sol-NLS-fourier}.
Then $A_n$ satisfy the system of ODEs
\begin{equation}
    \label{eq:A-ODE}
    \frac{\diff}{\diff t} A_n(t) = - i\epsilon \gamma \sum_{n_1,n_2,n_3}  \braket{\psi_n \psi_{n_1} \psi_{n_2} \psi_{n_3}} (A_{n_1} \overline{A_{n_2}} A_{n_3})(t)  e^{is \Omega(\lambda_n^L,\lambda_{n_1}^L,\lambda_{n_2}^L,\lambda_{n_3}^L)},
\end{equation}
where $\gamma = L^{-d} \vol(M)^{-1}$, $\langle \cdot \rangle = \langle \cdot \rangle_M$, and the resonance modulus $\Omega$ is given by
\begin{equation*}
    \Omega(\lambda,\lambda_1,\lambda_2,\lambda_3)
    = \lambda - \lambda_1 + \lambda_2 - \lambda_3.
\end{equation*}

\subsection{Wick ordering}

Denote $c_n = \varphi(\lambda_n^L) e^{i\theta_n}$ for simplicity.
We have
\begin{equation*}
    \mathfrak{m} \coloneq \|u\|_{L^2(M_L)}^2
    = \sum_{n \ge 0} |A_n(t)|^2 = \sum_{n \ge 0} |\varphi(\lambda_n^L)|^2
    \sim \zeta \int \omega^{d/2-1} |\varphi({\omega})|^2 \diff \omega,
    \quad \zeta = \frac{dL^d}{2 C_M (4\pi)^{d/2}}.
\end{equation*}
This is a consequence of the conservation of mass and the following estimate which holds true for all $\chi \in C_c^\infty(\mathbb{R}_+)$ and is proved using the Weyl law~\eqref{eq::weyl-law} and the convergence of Riemann sums:
\begin{equation}
    \label{eq::sum-to-integral}
    \sum_n \chi(\lambda_n^L)
    \sim \sum_n \chi\Bigl(\frac{ 4\pi (C_M n)^{2/d}}{L^2}\Bigr) \\
    \sim \zeta \int \omega^{d/2-1} \chi({\omega}) \diff \omega.
\end{equation}

We perform the Wick ordering by introducing $B_k(t) = A_k(t) e^{2i\varepsilon \gamma \mathfrak{m} t}$.
By~\eqref{eq:A-ODE}, one has
\begin{equation*}
    \frac{\diff}{\diff t} B_n(t) =
    -i\epsilon \gamma \Bigl( \sum_{n_1,n_2,n_3}^\times \braket{\psi_n \psi_{n_1} \psi_{n_2} \psi_{n_3}} (B_{n_1} \overline{B_{n_2}} B_{n_3})(t)  e^{it \Omega(\lambda_n^L,\lambda_{n_1}^L,\lambda_{n_2}^L,\lambda_{n_3}^L)} + R_n(t) B_n(t) \Bigr),
\end{equation*}
where $\sum_{n_1,n_2,n_3}^\times$ sums over all non-degenerate $(n_1,n_2,n_3)$, i.e., those with $(n,n_2) \ne (n_1,n_3)$ and $(n,n_2) \ne (n_3,n_1)$, and $R_n(t)$ encodes the contributions from the degenerate terms.
By~\eqref{eq::Berry-multi-eigen-z=0},
\begin{equation}
    \label{eq:est-Wick-ordering}
    R_n(t) = \sum_\ell \braket{\psi_n^2 \psi_\ell^2} |B_\ell(t)|^2 - 2 \mathfrak{m} \sim 0.
\end{equation}

\subsection{Perturbative expansion}

Using the Picard iteration method and Duhamel's principle, we obtain an asymptotic expansion $B_n \sim \sum_{j \ge 0} (-i\epsilon \gamma)^j B_n^j$, where the first few iterates are given by
\begin{align*}
    B_n^0(t) & = c_n; \\
    B_n^1(t) & \sim \sum_{n_1,n_2,n_3}^\times \braket{\psi_n \psi_{n_1} \psi_{n_2} \psi_{n_3}} (c_{n_1} \overline{c_{n_2}} c_{n_3}) \int_0^t e^{is \Omega} \diff s;                                                 \\
    B_k^2(t) & \sim \sum_{n_1,n_2,n_3}^\times \sum_{\ell_1,\ell_2,\ell_3}^\times \braket{\psi_n \psi_{n_1} \psi_{n_2} \psi_{n_3}}\braket{\psi_{n_1} \psi_{\ell_1} \psi_{\ell_2} \psi_{\ell_3}}
    (c_{\ell_1} \overline{c_{\ell_2}} c_{\ell_3}
    \overline{c_{n_2}} c_{n_3})
    \int_0^t \int_0^s e^{is\Omega} e^{is' \Omega_1} \diff s' \diff s \\
             & \quad
    + \sum_{n_1,n_2,n_3}^\times \sum_{\ell_1,\ell_2,\ell_3}^\times \braket{\psi_n \psi_{n_1} \psi_{n_2} \psi_{n_3}} \braket{\psi_{n_3} \psi_{\ell_1} \psi_{\ell_2} \psi_{\ell_3}}
    (c_{\ell_1} \overline{c_{\ell_2}} c_{\ell_3}
    \overline{c_{n_2}} c_{n_3})
    \int_0^t \int_0^s e^{is\Omega} e^{is' \Omega_3} \diff s' \diff s  \\
             & \quad
    - \sum_{n_1,n_2,n_3}^\times \sum_{\ell_1,\ell_2,\ell_3}^\times \braket{\psi_n \psi_{n_1} \psi_{n_2} \psi_{n_3}} \braket{\psi_{n_2} \psi_{\ell_1} \psi_{\ell_2} \psi_{\ell_3}}
    (c_{n_1} \overline{c_{\ell_1}} c_{\ell_2} \overline{c_{\ell_3}} c_{n_3}) \int_0^t \int_0^s e^{is\Omega} e^{-is' \Omega_2} \diff s' \diff s.
\end{align*}
where we denote $\Omega = \Omega(\lambda_n^L,\lambda_{n_1}^L,\lambda_{n_2}^L,\lambda_{n_3}^L)$, and $\Omega_j = \Omega(\lambda_{n_j}^L,\lambda_{\ell_1}^L,\lambda_{\ell_2}^L,\lambda_{\ell_3}^L)$, $j = 1,2,3$, for simplicity.

\subsection{Energy density}

We look for the asymptotics of the energy density $\mathbb{E}|A_n(t)|^2 = \mathbb{E}|B_n(t)|^2$. By the asymptotic expansion, we write
\begin{align*}
    |B_n|^2
     & \sim \sum_j (\epsilon\gamma)^j \sum_{p+q = j} (-i)^{p} i^{q} B_n^p \overline{B_n^q}                                                                                                                     \\
     & \sim |B_n^0|^2 -i\epsilon\gamma (B_n^0 \overline{B_n^1} + B_n^1 \overline{B_n^0}) -(\epsilon\gamma)^2 (B_n^1 \overline{B_n^1} - B_n^0 \overline{B_n^2} - B_n^2 \overline{B_n^0}) + \cdots
\end{align*}

Clearly $\mathbb{E}|B_n^0(t)|^2 = \mathbb{E}|c_n|^2 = |\varphi(\lambda_n^L)|^2$.
Next, note that
\begin{equation*}
    \mathbb{E}[\overline{c_n} c_{n_1} \overline{c_{n_2}} c_{n_3}] = \bm{1}_{(n,n_2) \in \{(n_1,n_3),(n_3,n_1)\}} |\varphi(\lambda_n^L)|^2 |\varphi(\lambda_{n_2}^L)|^2,
\end{equation*}
which clearly vanishes for all $n_1,n_2,n_3$ appearing in the summation $\sum^\times_{n_1,n_2,n_3}$.
Therefore
\begin{equation*}
    \mathbb{E}[B_n^0 \overline{B_n^1} + B_n^1 \overline{B_n^0}](t)
    \sim 2 \Re \Bigl(\sum_{n_1,n_2,n_3}^\times \braket{\psi_n \psi_{n_1} \psi_{n_2} \psi_{n_3}} \mathbb{E}[\overline{c_n} c_{n_1} \overline{c_{n_2}} c_{n_3}] \int_0^t  e^{is \Omega} \diff s \Bigr) = 0.
\end{equation*}
Again, by the Wick ordering, the estimate~\eqref{eq:est-interaction-real}, and the scaling relation~\eqref{eq::scaling-I}, we obtain
\begin{align*}
    \mathbb{E}|B_n^1|^2(t)
    & =\sum_{n_1,n_2,n_3}^\times \sum_{\ell_1,\ell_2,\ell_3}^\times
    \braket{\psi_n \psi_{n_1} \psi_{n_2} \psi_{n_3}}
    \braket{\psi_n \psi_{\ell_1} \psi_{\ell_2} \psi_{\ell_3}}
    \mathbb{E}[c_{n_1} \overline{c_{n_2}} c_{n_3}
        \overline{c_{\ell_1}} c_{\ell_2} \overline{c_{\ell_3}}]
    \int_0^t  e^{is \Omega} \diff s
    \int_0^t e^{-is' \Omega_1} \diff s'\\
    & \sim 2 \sum_{n_1,n_2,n_3}^\times \braket{\psi_n \psi_{n_1} \psi_{n_2} \psi_{n_3}}^2 |\varphi(\lambda_{n_1}^L)|^2 |\varphi(\lambda_{n_2}^L)|^2 |\varphi(\lambda_{n_3}^L)|^2 \Bigl| \int_0^t e^{is \Omega} \diff s \Bigr|^2\\
    & \sim 2 \vol(M)^{-1} \sum_{n_1,n_2,n_3}^\times \mathcal{I}(\lambda_{n},\lambda_{n_1},\lambda_{n_2},\lambda_{n_3}) |\varphi(\lambda_{n_1}^L)|^2 |\varphi(\lambda_{n_2}^L)|^2 |\varphi(\lambda_{n_3}^L)|^2 \Bigl| \int_0^t e^{is \Omega} \diff s \Bigr|^2\\
     & \sim 2 \gamma \sum_{n_1,n_2,n_3}^\times \mathcal{I}(\lambda_{n}^L,\lambda_{n_1}^L,\lambda_{n_2}^L,\lambda_{n_3}^L) |\varphi(\lambda_{n_1}^L)|^2 |\varphi(\lambda_{n_2}^L)|^2 |\varphi(\lambda_{n_3}^L)|^2 \Bigl| \int_0^t e^{is \Omega} \diff s \Bigr|^2.
\end{align*}
Similarly, identifying $n=n_0$, we obtain that
\begin{equation*}
    \mathbb{E} [B_n^0 \overline{B_n^2} + B_n^2 \overline{B_n^0}](t)
    \sim -2 \gamma \sum_{n_1,n_2,n_3}^\times \Bigl(\sum_{1\le \ell\le 3}(-1)^\ell \prod_{0\le j\le 3}^{j\ne \ell} |\varphi(\lambda_{n_j}^L)|^2\Bigr) \mathcal{I}(\lambda_n^L,\lambda_{n_1}^L,\lambda_{n_2}^L,\lambda_{n_3}^L) \Bigl| \int_0^t e^{is \Omega} \diff s \Bigr|^2.
\end{equation*}
Therefore, neglecting the lower order terms and using $\int_0^t e^{is \Omega} \diff s = e^{it\Omega/2} \frac{\sin(t\Omega/2)}{\Omega/2}$, we have
\begin{equation*}
    \mathbb{E}|B_n|^2(t)
     \sim |\varphi(\lambda_n^L)|^2 + 2 \epsilon^2 \gamma^3 \sum_{n_1,n_2,n_3} F(\lambda_n^L,\lambda_{n_1}^L,\lambda_{n_2}^L,\lambda_{n_3}^L) \mathcal{I}(\lambda_{n}^L,\lambda_{n_1}^L,\lambda_{n_2}^L,\lambda_{n_3}^L) \frac{\sin(t\Omega/2)^2}{(\Omega/2)^2},
\end{equation*}
where $F(\lambda_n^L,\lambda_{n_1}^L,\lambda_{n_2}^L,\lambda_{n_3}^L) = \sum_{0\le \ell\le 3}(-1)^\ell \prod_{j\ne \ell} |\varphi(\lambda_{n_j}^L)|^2.$

\subsection{Kinetic limit}

\label{sec::kinetic-limit}

Pass to the integral using~\eqref{eq::sum-to-integral}, write $\omega = \lambda_n^L$, and with an abuse of notation denote $\Omega = \Omega(\omega,\omega_1,\omega_2,\omega_3)$.
The approximation $\frac{\sin(t\Omega/2)^2}{(\Omega/2)^2} \sim 2\pi t \bm{\delta}_{\mathbb{R}}(\Omega)$ for $t \gg 1$ then yields
\begin{align*}
    \mathbb{E}|B_n|^2(t) - |\varphi(\omega)|^2
        & \sim 2\epsilon^2\gamma^3 \zeta^3 \iiint_{(\mathbb{R}_+)^3} F(\omega,\omega_1,\omega_2,\omega_3) \mathcal{I}(\omega,\omega_1,\omega_2,\omega_3) \frac{\sin(t\Omega/2)^2}{(\Omega/2)^2} \diff \omega_1 \diff \omega_2 \diff \omega_3 \\
        & \sim 4\pi t\epsilon^2\gamma^3 \zeta^3 \iiint_{(\mathbb{R}_+)^3} F(\omega,\omega_1,\omega_2,\omega_3) \mathcal{I}(\omega,\omega_1,\omega_2,\omega_3) \bm{\delta}_{\mathbb{R}}(\Omega)\diff \omega_1 \diff \omega_2 \diff \omega_3.
\end{align*}
Our choice of the kinetic time order~\eqref{eq::Tkin} is then justified by the computation:
\begin{equation*}
    4\pi t\epsilon^2\gamma^3 \zeta^3
    = \frac{\pi^2}{2} \times \biggl(\frac{s(d)}{(2\pi)^d}\biggr)^3 \times \frac{\varepsilon^2 t}{\pi}.
\end{equation*}

\bibliographystyle{abbrv}
\bibliography{references}

\end{document}